\documentclass{amsart}

\usepackage{amsmath}
\usepackage{amsthm}
\usepackage{amsfonts}
\usepackage{amssymb}

\usepackage{graphics}
\usepackage{epsfig}
\usepackage{color}
\usepackage{verbatim}
\newtheorem{theorem}{Theorem}
\newtheorem{lemma}[theorem]{Lemma}

\newtheorem{corollary}[theorem]{Corollary}

\theoremstyle{remark}
\newtheorem*{remark}{Remark}
\newtheorem*{question*}{Question}

\date{20 May 2019}
\title[A supersolutions perspective on hypercontractivity]{A supersolutions perspective on hypercontractivity}
\begin{document}

\author[Y. Aoki]{Yosuke Aoki}
\address{Department of Mathematics, Graduate School of Science and Engineering,
Saitama University, Saitama 338-8570, Japan.}
\email{y.aoki.393@ms.saitama-u.ac.jp (Aoki), \newline nealbez@mail.saitama-u.ac.jp (Bez), \newline machihara@rimath.saitama-u.ac.jp (Machihara), \newline k.matsuura.490@ms.saitama-u.ac.jp (Matsuura), \newline s.shiraki.446@ms.saitama-u.ac.jp (Shiraki)}

\author[J. Bennett]{Jonathan Bennett}
\address{School of Mathematics, The Watson Building, University of Birmingham, Edgbaston, Birmingham, B15 2TT, England.}
\email{J.Bennett@bham.ac.uk (Bennett)}

\author[N. Bez]{Neal Bez}

\author[S. Machihara]{Shuji Machihara}

\author[K. Matsuura]{Kosuke Matsuura}

\author[S. Shiraki]{Shobu Shiraki}

\thanks{The first, fifth and sixth authors were supported by JSPS Grant-in-Aid for Young Scientists A [grant number 16H05995], the second author was partially supported by ERC grant 307617, the third author was supported by JSPS Grant-in-Aid for Young Scientists A [grant number 16H05995] and JSPS Grant-in-Aid for Scientific Research B [grant number 19H01796], and the fourth author was supported by JSPS Grant-in-Aid for Scientific Research C [grant number 16K05191].}

\begin{abstract}
The purpose of this article is to expose an algebraic closure property of supersolutions to certain diffusion equations. This closure property quickly gives rise to a monotone quantity which generates a hypercontractivity inequality. Our abstract argument applies to a general Markov semigroup whose generator is  a   diffusion and satisfies a curvature condition.
\end{abstract}
\maketitle

\section{Introduction}

We begin in the concrete setting of the Ornstein--Uhlenbeck semigroup $(e^{sL})_{s \geq 0}$ defined by
\begin{equation} \label{e:OUsemi}
e^{sL}f(x) = \int_{\mathbb{R}^d} f(e^{-s}x + (1-e^{-2s})^{1/2}y)\, \mathrm{d}\gamma(y),
\end{equation}
where $d \geq 1$ and $\gamma$ is the standard gaussian probability measure on $\mathbb{R}^d$ given by
\[
\mathrm{d}\gamma(y) = e^{-\frac{1}{2}|y|^2}  \frac{\mathrm{d}y}{(2\pi)^{d/2}}.
\]
The generator is the so-called number operator given by $Lf(x) = \Delta f(x) - x \cdot \nabla f(x)$.
\begin{theorem} \label{t:closure_OU}
Let $\infty > q > p > 1$ and $s > 0$ be given by $e^{2s} = \frac{q-1}{p-1}$. Suppose $u : (0,\infty) \times \mathbb{R}^d \to (0,\infty)$ is such that $u(t,\cdot)^{1/p}$, $\partial_t(u(t,\cdot)^{1/p})$, $\nabla(u(t,\cdot)^{1/p})$, $u(t,\cdot)^{-1/p}|\nabla(u(t,\cdot)^{1/p})|^2$ and $\Delta(u(t,\cdot)^{1/p})$ are of polynomial growth locally uniformly in time $t > 0$, and satisfies
\[
\partial_t u \geq L u.
\]
Let $\widetilde{u} : (0,\infty) \times \mathbb{R}^d \to (0,\infty)$ be given by
\begin{equation} \label{e:utilde}
\widetilde{u}(t,x)^{1/q} = e^{sL}(u(t,\cdot)^{1/p})(x).
\end{equation}
Then $\widetilde{u}(t,\cdot)^{1/q}$, $\partial_t(\widetilde{u}(t,\cdot)^{1/q})$, $\nabla(\widetilde{u}(t,\cdot)^{1/q})$, $\widetilde{u}(t,\cdot)^{-1/q}|\nabla(\widetilde{u}(t,\cdot)^{1/q})|^2$ and $\Delta(\widetilde{u}(t,\cdot)^{1/q})$ are of polynomial growth locally uniformly in time $t > 0$, and 
\[
\partial_t \widetilde{u} \geq L \widetilde{u}.
\]
\end{theorem}
The main feature of this theorem is that supersolutions of the linear diffusion equation governed by $L$ are preserved under the transformation $u \mapsto \widetilde{u}$. Although the regularity conditions imposed on $u$ are of a more technical nature, some care has been exercised to ensure that they are strong enough for the relevant terms in the statement of the theorem and its proof to be rigorously defined, and weak enough so that the regularity conditions are themselves preserved under the transformation $u \mapsto \widetilde{u}$. In this regard, our theorem is compatible with the perspective taken in \cite{BB_Crelle}, where a framework for generating monotone quantities for nonnegative solutions of linear heat equations was developed based on algebraic closure properties of supersolutions.

The monotone quantity which most immediately arises from Theorem \ref{t:closure_OU} is contained in the following.
\begin{corollary} \label{c:mono_OU}
Suppose $u$ satisfies $\partial_t u = L u$ with initial data a bounded and compactly supported nonnegative function on $\mathbb{R}^d$.  Let $Q : (0,\infty) \to (0,\infty)$ be given by
\[
Q(t) = \int_{\mathbb{R}^d} \left(e^{sL}(u(t,\cdot)^{1/p})\right)^q(x) \, \mathrm{d}\gamma(x),
\] 
where $\infty > q > p > 1$ and $e^{2s} = \frac{q-1}{p-1}$. Then $Q$ is nondecreasing on $(0,\infty)$.
\end{corollary}
Taking $\widetilde{u}$ as in \eqref{e:utilde}, we have $Q(t) = \int \widetilde{u}(t,\cdot) \, \mathrm{d}\gamma$ and by passing the time derivative through the integral, we may quickly obtain Corollary \ref{c:mono_OU} from Theorem \ref{t:closure_OU}. In turn, the monotonicity of $Q$ generates the well-known hypercontractivity inequality enjoyed by the Ornstein--Uhlenbeck semigroup. Indeed, taking $u$ to satisfy $\partial_t u = L u$ with initial data $f^p$, where $f$ is a bounded and compactly supported nonnegative function on $\mathbb{R}^d$, the dominated convergence theorem implies that
\[
\lim_{t \to 0} Q(t) = \| e^{sL} f \|_{L^q(\gamma)}^q
\]
and
\[
\lim_{t \to \infty} Q(t) = \|f\|_{L^p(\gamma)}^q.
\]
In this manner, Theorem \ref{t:closure_OU} quickly yields Nelson's famous hypercontractivity inequality 
\begin{equation} \label{e:hyper_OU}
\| e^{sL} f \|_{L^q(\gamma)} \leq \|f\|_{L^p(\gamma)}
\end{equation}
whenever $q > p > 1$ and $s > 0$ are such that
\begin{equation} \label{e:Nelsontime}
e^{2s} = \frac{q-1}{p-1}.
\end{equation}
For such $p$ and $q$, the time $s$ given by \eqref{e:Nelsontime} is critical since the operator $e^{sL}$ is unbounded from $L^p(\gamma)$ to $L^q(\gamma)$ for any smaller value of $s$. Nelson \cite{Nelson} first derived the inequality in \eqref{e:hyper_OU} in his work on quantum field theory and later Gross \cite{Gross_AJM} established an equivalence with certain log-Sobolev inequalities. For further historical details and wider perspectives on the role played by this inequality and its generalisations, we refer the reader to \cite{DGS} and \cite{Gross_survey}. We also note that different proofs of the hypercontractivity inequality based on monotone quantities may already be found in, for example, work of Hu \cite{Hu} and Ledoux \cite{Ledoux}.

The key argument on which our proof of Theorem \ref{t:closure_OU} is based may be applied, at least at a formal level, in the significantly more general setting of Markov semigroups. We consider this to be our main contribution in the current work and in order to expose this in the clearest manner, in the forthcoming Section \ref{section:abstract}, we first present the abstract argument. The argument in Section \ref{section:abstract} is used in Section \ref{section:OU} to prove Theorem \ref{t:closure_OU} and Corollary \ref{c:mono_OU}. Finally, in Section \ref{section:further} make some further remarks concerning closure properties associated with reverse hypercontractivity.

\section{The abstract argument} \label{section:abstract}

\subsection{Some preparation} The underlying setting is a $\sigma$-finite measure space $(E,\mathcal{E},\mu)$ and a Markov semigroup $(P_s)_{s \geq 0}$ given by 
\begin{equation} \label{e:repformula}
P_s f(x)=\int_E f(y) \, \mathrm{d}\nu_{s,x}(y)
\end{equation}
for $x \in E$, where $\nu_{s,x}$ is a nonnegative probability measure (``transition kernel"). Associated to the semigroup is the Markov generator $L$ and we assume that the underlying measure $\mu$ is invariant with respect to $L$. All operations are assumed to be well-defined on an appropriate algebra of functions on $E$. Associated to $L$ are the \emph{carr\'e du champ operator} and the \emph{curvature operator} given by
\[
\Gamma(f,g)=\frac{1}{2}(L(fg)-fLg-gLf)
\]
and
\[
\Gamma_2(f,g)=\frac{1}{2}(L\Gamma(f,g)-\Gamma(f,Lg)-\Gamma(g,Lf)),
\]
respectively. We set $\Gamma(f)=\Gamma(f,f)$ and similarly for $\Gamma_2$. 

We say $L$ is a \emph{diffusion} if, for all $C^\infty$ functions $\psi$ on $\mathbb{R}^n$, we have
\begin{align}\label{e:diffusion}
L\psi(f) = \sum_{j=1}^n \partial_j\psi(f) Lf_j +  \sum_{j,k = 1}^n \partial_{j,k}^2\psi(f) \Gamma(f_j,f_k), 
\end{align}
where $f = (f_1,\ldots,f_n)$, and of \emph{curvature $c \in \mathbb{R}$} if 
\begin{equation} \label{e:curvature}
\Gamma_2(f) \geq c\Gamma(f).
\end{equation}
These properties suffice\footnote{In fact, if the diffusion property holds, then the curvature condition is \emph{equivalent} to \eqref{e:convexity}.} for the following key property to hold. 
\begin{lemma}\label{l:convexity}
If $L$ is a diffusion and of curvature $c$, then
\begin{equation} \label{e:convexity}
\sqrt{\Gamma (P_s f)} \leq e^{-cs}P_s[\sqrt{\Gamma(f)}]
\end{equation}
for all $s\geq0$.
\end{lemma}
The above gradient bound is due to Bakry \cite{Bakry}. We also refer the reader to \cite{BGL} and \cite{Ledoux_Survey} for further details regarding the abstract setting we are working in.
\begin{remark}
In the case of the Ornstein--Uhlenbeck semigroup $P_s f = e^{sL}f$ where $L = \Delta - x \cdot \nabla$, a simple change of variables shows that \eqref{e:repformula} holds with
\[
\mathrm{d}\nu_{s,x}(y)= \exp\bigg(-\frac{|\rho x - y|^2}{2(1-\rho^2)} \bigg)  \frac{\mathrm{d}y}{[2\pi(1-\rho^2)]^{d/2}},
\]
where $\rho = e^{-s}$. Also, direct computations reveal that $\Gamma(f) = |\nabla f|^2$ and $\Gamma_2(f) = |D^2f|^2 + |\nabla f|^2$; thus, \eqref{e:diffusion} and \eqref{e:curvature} with $c=1$ hold. In this special case, the explicit formula \eqref{e:OUsemi} immediately implies the key estimate \eqref{e:convexity}.
\end{remark}

\subsection{The closure property}
Suppose $\infty > q > p > 1$ and $s$ is defined by $e^{2cs} = \frac{q-1}{p-1}$, where $c \in \mathbb{R}$ is the curvature constant in \eqref{e:curvature}. Let $\widetilde{u} : (0,\infty) \times \mathbb{R}^d \to (0,\infty)$ be given by
\[
\widetilde{u}(t,x)^{1/q} = P_s[u(t,\cdot)^{1/p}](x)
\]
for $(t,x) \in (0,\infty) \times \mathbb{R}^d$. We shall prove here that $\partial_t \widetilde{u} \geq L\widetilde{u}$ whenever $u$ satisfies $\partial_t u \geq Lu$; for simplicity of the exposition, the following argument is based on certain formal considerations. For instance, we shall make multiple use of the identity
\begin{equation} \label{e:diffusion_powers}
L(f^ \lambda)= \lambda f^{\lambda-1}Lf + \lambda(\lambda-1)f^{\lambda-2}\Gamma(f)
\end{equation}
for $\lambda > 0$. Observe that \eqref{e:diffusion_powers} formally follows from the diffusion property by taking $\psi(f) = f^\lambda$ (or in a rigorous sense in the case of the Ornstein--Uhlenbeck semigroup via direct computations).

Proceeding via the representation formula \eqref{e:repformula} and formally passing the time derivative through the integral, we have
\begin{equation} \label{e:formal1}
\partial_t\widetilde{u}(t,x) = \frac{q}{p} \widetilde{u}(t,x)^{1-1/q} P_s[u(t,\cdot)^{1/p-1}\partial_tu(t,\cdot)](x).
\end{equation}
Since $u$ is a supersolution and by use of \eqref{e:diffusion_powers} we get
\begin{align*}
\partial_t\widetilde{u}(x) & \geq \frac{q}{p} \widetilde{u}(x)^{1-1/q} P_s[u^{1/p-1}Lu](x) \\
&=q \widetilde{u}(x)^{1 - 1/q} P_s[Lu^{1/p}](x) +\frac{q}{pp'} \widetilde{u}(x)^{1 - 1/q} P_s[u^{1/p-2}\Gamma(u)](x).
\end{align*}
Here we have dropped the dependence on the $t$ variable since all operators are now acting in the spatial variable.

On the other hand, by a further application of \eqref{e:diffusion_powers},
\begin{equation*}
L\widetilde{u}(x) =q \widetilde{u}(x)^{1 - 1/q} LP_s[u^{1/p}](x)+q(q-1)\widetilde{u}(x)^{1 - 2/q}\Gamma(P_s[u^{1/p}])(x) 
\end{equation*}
and, using that $P_s$ and $L$ formally commute, we thus have
\begin{equation*}
\frac{1}{q}\widetilde{u}(x)^{2/q - 1}[\partial_t\widetilde{u} - L\widetilde{u}](x) \geq \frac{1}{pp'} \widetilde{u}(x)^{1/q} P_s[u^{1/p-2}\Gamma(u)](x) - (q-1)\Gamma(P_s[u^{1/p}])(x). 
\end{equation*}
However, by an application of Lemma \ref{l:convexity} followed by the Cauchy--Schwarz inequality we have
\begin{align*}
(q-1)\Gamma(P_s[u^{1/p}])(x) \leq (q-1)e^{-2cs} P_s[u^{1/p}](x) P_s[u^{-1/p}\Gamma(u^{1/p})](x).
\end{align*}
Applying the identity\footnote{The identity $\Gamma(\psi(u))= \psi'(u)^2\Gamma(u)$ holds for smooth $\psi$ as a result of the diffusion property \eqref{e:diffusion} and thus \eqref{e:property_gamma} holds in a formal sense by taking $\psi(u) = u^{1/p}$. In the case of the Ornstein--Uhlenbeck semigroup, \eqref{e:property_gamma} may be rigorously verified by direct calculations.}
\begin{equation} \label{e:property_gamma}
\Gamma(u^{1/p}) = \frac{1}{p^2} u^{2/p - 2} \Gamma(u)
\end{equation}
and using the relation $e^{2cs} = \frac{q-1}{p-1}$, it is clear from the above argument that $\partial_t\widetilde{u} \geq L\widetilde{u}$.

\section{The Ornstein--Uhlenbeck semigroup} \label{section:OU}

In this section, we write $P_sf = e^{sL}f$, where $L = \Delta - x \cdot \nabla$, and $B_s(x,y) = e^{-s}x + (1-e^{-2s})^{1/2}y$.

\subsection{Proof of Theorem \ref{t:closure_OU}}

To begin, we observe that $\widetilde{u}(t,x)$ is well-defined in a pointwise sense since our assumptions on $u$ mean that $u(t,\cdot)^{1/p}$ is of polynomial growth for each fixed time.

In order to prove $\partial_t\widetilde{u} \geq L\widetilde{u}$, we run the argument in the previous section with $c=1$. Rigorous justification of \eqref{e:formal1}, at which point we passed the time derivative through the integral appearing in \eqref{e:repformula}, is made using the fact that $\partial_t (u^{1/p})$ is of polynomial growth locally uniformly in time. Another formal step in the argument is the commutativity property 
\[
P_s[L(u^{1/p})](x) = LP_s[u^{1/p}](x).
\]
Since $L = \Delta - x \cdot \nabla$, we may rigorously justify this since $\nabla(u^{1/p})$ and $\Delta(u^{1/p})$ are of polynomial growth locally uniformly in time $t > 0$. Finally, we note that the term $P_s[u^{1/p-2}\Gamma(u)](x)$ is well-defined in a pointwise sense thanks to the assumption that $u^{-1/p}|\nabla(u^{1/p})|^2$ is of polynomial growth locally uniformly in time. This completes the verification of the formal steps in the argument in the previous section.

It remains to check that the regularity conditions imposed on $u$ result in $\widetilde{u}$ satisfying analogous regularity properties. Our assumption on $u$ means that, for a fixed $t > 0$, there is a natural number $N$ and compact interval $I \subset (0,\infty)$ containing $t$ such that $\sup_{t' \in I} u(t',x)^{1/p} \lesssim_t \langle x \rangle^N$ for all $x \in \mathbb{R}^d$. Here, we are using the Japanese bracket notation $\langle x \rangle = (1 + |x|^2)^{1/2}$. Thus, clearly we have
\begin{equation} \label{e:techincal_u}
\widetilde{u}(t',x)^{1/q} = \int u^{1/p}(t',B_s(x,y)) \, \mathrm{d}\gamma(y) \lesssim_t \int \langle B_s(x,y) \rangle^{N} \, \mathrm{d}\gamma(y) \lesssim_{s,t} \langle x \rangle^{N}
\end{equation}
for each $t' \in I$ and $x \in \mathbb{R}^d$, and it follows that $\widetilde{u}^{1/q}$ is of polynomial growth locally uniformly in time.

For $\partial_t(\widetilde{u}(t,\cdot)^{1/q})$, the assumption that $\partial_t (u^{1/p})$ is of polynomial growth locally uniformly in time means, by a routine application of the dominated convergence theorem,
\[
\partial_t(\widetilde{u}^{1/q})(t',x) = \int \partial_t(u^{1/p})(t',B_s(x,y)) \, \mathrm{d}\gamma(y) 
\]
for all $t'$ in an appropriate compact interval, and now estimating in a similar manner to \eqref{e:techincal_u} reveals that $\partial_t(\widetilde{u}(t,\cdot)^{1/q})$ of polynomial growth locally uniformly in time. By similar considerations, the same conclusion also holds for $\nabla(\widetilde{u}^{1/q})$ and $\Delta(\widetilde{u}^{1/q})$. Finally, by the Cauchy--Schwarz inequality
\begin{align*}
|\nabla(\widetilde{u}^{1/q})(t',x)|^2 & = \bigg| \int \nabla(u^{1/p})(t',B_s(x,y)) \, \mathrm{d}\gamma(y) \bigg|^2 \\
& \leq \widetilde{u}(t',x)^{1/q} \int \frac{|\nabla(u^{1/p}(t',B_s(x,y))|^2}{u^{1/p}(t',B_s(x,y))} \, \mathrm{d}\gamma(y)
\end{align*}
and the fact that $u^{-1/p}|\nabla(u^{1/p})|^2$ is of polynomial growth locally uniformly in time can be easily seen to induce the same property for $\widetilde{u}^{-1/q} |\nabla(\widetilde{u}^{1/q})|^2$.

\subsection{Proof of Corollary \ref{c:mono_OU}}
Suppose that $f$ is bounded and nonnegative function with support inside $\{x \in \mathbb{R}^d : |x| \leq R\}$, and let $u(t,x) = P_t[f^p](x)$. By Theorem \ref{t:closure_OU}, it suffices to show that $u(t,\cdot)^{1/p}$, $\partial_t(u(t,\cdot)^{1/p})$, $\nabla(u(t,\cdot)^{1/p})$, $u(t,\cdot)^{-1/p}|\nabla(u(t,\cdot)^{1/p})|^2$ and $\Delta(u(t,\cdot)^{1/p})$ are of polynomial growth locally uniformly in time $t > 0$. Indeed, if this is the case, then $\partial_t \widetilde{u} \geq L \widetilde{u}$ where $\widetilde{u}(t,x)^{1/q} = P_s(u(t,\cdot)^{1/p})(x)$ and therefore
\[
Q'(t) = \frac{1}{q} Q(t)^{1 - 1/q} \int \partial_t \widetilde{u}(t,x) \, \mathrm{d}\gamma(x) \geq \frac{1}{q} Q(t)^{1 - 1/q} \int L \widetilde{u}(t,x) \, \mathrm{d}\gamma(x) = 0.
\]
Note that we may use the dominated convergence theorem to justify the interchange of the time derivative and the integral in the above argument. Indeed, we know from Theorem \ref{t:closure_OU} that both $\widetilde{u}(t,\cdot)^{1/q}$ and $\partial_t(\widetilde{u}(t,\cdot)^{1/q})$ are of polynomial growth locally uniformly in $t > 0$. By writing $\partial_t \widetilde{u} = q \widetilde{u}^{1-1/q} \partial_t(\widetilde{u}^{1/q})$ and recalling that $q > 1$, we see that the same property also holds for $\widetilde{u}$, and this is sufficient to justify the interchange of time derivative and integral.

It remains verify the regularity claimed hypotheses for $u$. We first note that $u(t,x)^{1/p} \leq \|f\|_\infty$ obviously follows from \eqref{e:OUsemi}. For $\partial_t(u^{1/p})$, we shall make use of the representation formula
\[
u(t,x) = C(t) \int_{\mathbb{R}^d} f(y)^p  \exp\bigg(-\frac{|\rho(t) x - y|^2}{2(1-\rho(t)^2)} \bigg)  \,\mathrm{d}y,
\]
where $C(t) = [2\pi(1-\rho(t)^2)]^{-d/2}$ and $\rho(t) = e^{-t}$. Using the assumption on the support of $f$, it easily follows that
\[
|\partial_t u(t,x)| \lesssim_R (1-\rho(t)^2)^{-2} \langle x \rangle^2 u(t,x)
\]
and therefore, since $\partial_t(u^{1/p}) = \frac{1}{p} \frac{\partial_t u}{u} u^{1/p}$, we see that $\partial_t (u^{1/p})$ is of polynomial growth locally uniformly in $t > 0$. A similar argument reveals $|\nabla u(t,x)| \lesssim_R (1-\rho(t)^2)^{-1} \langle x \rangle u(t,x)$. From this we quickly obtain that $\nabla(u^{1/p})$ is of polynomial growth locally uniformly in $t > 0$ and, via the identity $u^{-1/p}|\nabla(u^{1/p})|^2 = \frac{1}{p} u^{1/p} \frac{|\nabla u|^2}{u^2}$, the same conclusion too for $u^{-1/p}|\nabla(u^{1/p})|^2$. Finally, similar considerations show that $\Delta(u^{1/p})$ is also of polynomial growth locally uniformly in $t > 0$.

\section{Further remarks} \label{section:further}

A further appealing feature of our abstract argument in Section \ref{section:abstract} is that it applies to exponents $p$ and $q$ in the setting of the reverse hypercontractivity inequality. In general, we let $\widetilde{u}$ be given by
\[
\widetilde{u}(t,x)=
\begin{cases}
P_s[u(t,\cdot)^{1/p}]^q(x)  &\text{if $p,q\not=0$},\\
P_s[e^{u(t,\cdot)}]^q(x)  &\text{if $p=0,q\not=0$},\\
\log{P_s[u(t,\cdot)^{1/p}]}(x)  &\text{if $p\not=0, q=0$}.
\end{cases}
\]
Then, at least in a formal sense, the following closure properties hold:
\begin{equation} \label{e:closure_general}
\begin{array}{lllll}
& \partial_t u \geq Lu \,\, \Rightarrow \,\, \partial_t \widetilde{u} \geq L\widetilde{u} \qquad \qquad & \text{for $1<p<q<\infty$, and for $-\infty<q < p < 0$,} \\
& \partial_t u \leq Lu \,\, \Rightarrow \,\, \partial_t \widetilde{u} \leq L\widetilde{u} \qquad \qquad & \text{for $0 \leq q < p < 1$,} \\
& \partial_t u \leq Lu \,\, \Rightarrow \,\, \partial_t \widetilde{u} \geq L\widetilde{u} \qquad \qquad & \text{for $-\infty< q < 0 \leq p < 1$.}
\end{array}
\end{equation}
As a result, in each of the above cases, one may obtain the monotonicity of
\begin{equation*}
Q(t) = \left\{ \begin{array}{lll}  (\int \widetilde{u}(t,\cdot) \,\mathrm{d}\mu)^{1/q} & \text{if $q\neq 0$} \\ \exp(\int \widetilde{u}(t,\cdot) \, \mathrm{d}\mu)  & \text{if $q = 0$}  \end{array}  \right.
\end{equation*}
for solutions $u$ of the diffusion equation $\partial_t u = Lu$ with nonnegative initial data.

The closure properties \eqref{e:closure_general} for $q < p < 1$ are associated with Borell's reverse form of the hypercontractivity inequality
\begin{equation} \label{e:lower}
\| P_s f \|_{L^q(\mu)} \geq \|f\|_{L^p(\mu)}
\end{equation}
for positive functions. For the Ornstein--Uhlenbeck semigroup, Borell \cite{Borell} observed that a unified approach to both forward and reverse hypercontractivity inequalities may be taken, whereby one first establishes a discrete ``Boolean hypercontractivity inequality" (the forward form independently due to Bonami \cite{Bonami} and Gross \cite{Gross_AJM}) and then applies the central limit theorem. Both the forward and reverse forms of the Boolean hypercontractivity inequality have also found numerous applications, notably in various fields of computer science; the reader is encouraged to look at \cite{Montanaro}, \cite{MORSS} and \cite{OD} for interesting examples and further references. We also remark that the theory of reverse hypercontractivity was significantly developed in recent work of Mossel--Oleszkiewicz--Sen \cite{MOS} along with a host of applications.

\end{document}